\documentclass[oneside,10pt]{article}
\usepackage[b5paper]{geometry}	

\usepackage{amsfonts,amsmath,latexsym,amssymb}
\usepackage{theorem}
\usepackage{mathrsfs,upref}
\usepackage{mathptmx}		

\usepackage{arXiv}	

\usepackage[dvipdf]{graphicx}

\newtheorem{statement}{Statement}

\newtheorem{lemma}{Lemma}
\newtheorem{corollary}{Corollary}

\theoremstyle{definition}

\newtheorem{example}{Example}

\newtheorem{remark}{Remark}

\begin{document}

\title[EMI]{On the Extension of the Erd\" os-Mordell Type Inequalities}


\author{B. Male\v sevi\' c, M. Petrovi\' c, M. Obradovi\' c, B. Popkonstantinovi\' c}

\address{} 

\date{04.04.2012.}



\thanks{Research is partially supported by the Ministry of Science and Education of the Republic of Serbia, Grant No. III~44006 and ON~174032}

\begin{abstract}
We discuss the extension of inequality $R_A \ge \frac{c}{a} r_b + \frac{b}{a}r_c$ to the plane of triangle $\triangle ABC$.
Based on the obtained extension, in regard to all three vertices of the triangle, we get the extension of Erd\" os-Mordell inequality,
and some inequalities of Erd\" os-Mordell type.
\end{abstract}

\maketitle

\section{Introduction}

Let triangle $\triangle ABC$ be given in Euclidean plane. Denote by $R_A, \, R_B$ and $R_C$ the distances from the arbitrary point $M$
in the interior of $\triangle ABC$ to the vertices $A, \, B$ and $C$ respectively, and denote by $r_a, \, r_b$ and $r_c$ the distances
from the point $M$ to the sides $BC, \, CA$ and $AB$ respectively (Figure 1).

\begin{center} 

\vspace*{20.0mm} \hspace*{-25.0mm} \includegraphics*[height=30.0mm,keepaspectratio=true,scale=0.65]{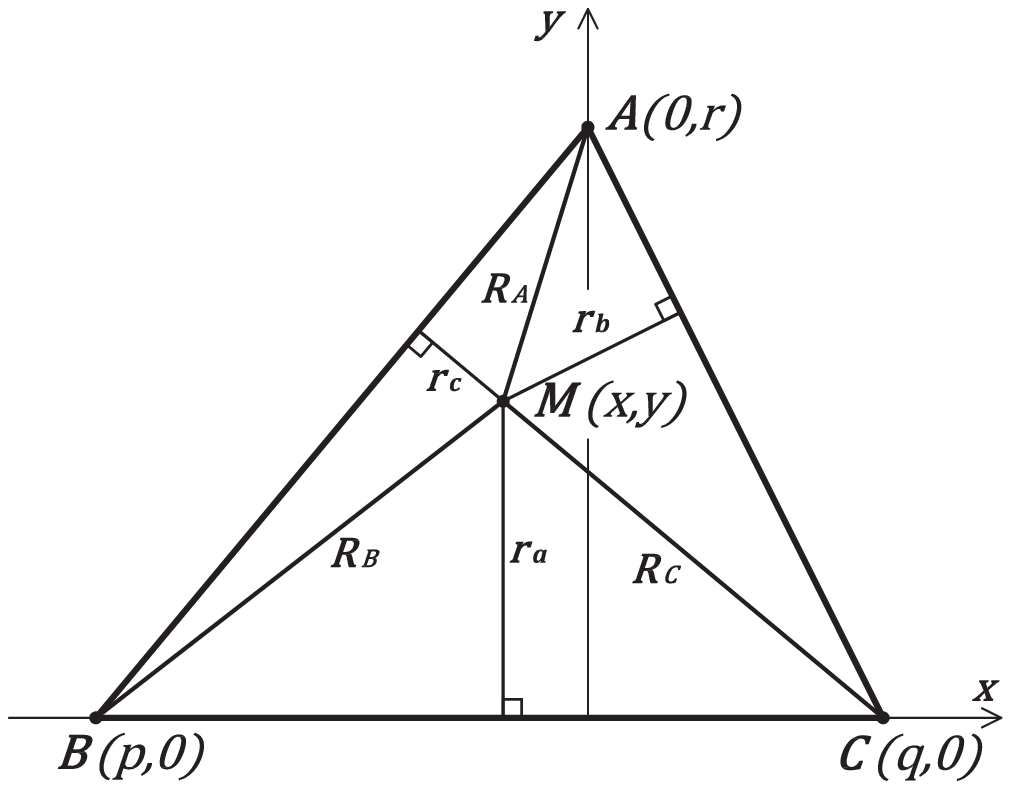}

\noindent
\textit{Figure 1: Erd\" os-Mordell inequality}

\end{center}

\smallskip
Then Erd\" os-Mordell inequality is true:
\begin{equation}
\label{GrindEQ__1_}
R_A + R_B + R_C
\ge
2\left( r_a + r_b + r_c \right)
\end{equation}
whereat equality holds if and only if triangle $ABC$ is equilateral and $M$ is its center. This inequality was conjectured by P. Erd\" os
as Amer. Math. Monthly Problem 3740 in 1935. \cite{[1]}, after his experimental conjecture in 1932. \cite{[13]}. It was proved by L.J. Mordell in 1935.
(in Hungarian, according to \cite{[13]}), and as the solution of the Problem 3740 in 1937. \cite{[2]}.

\break

\smallskip
Considering the Erd\" os-Mordell inequality \eqref{GrindEQ__1_} the goal of this research is to determine areas in the plane of the triangle,
where the following three inequalities are valid:
\begin{equation}
\label{GrindEQ__2_}
R_A \ge \frac{c}{a}r_b + \frac{b}{a}r_c
\end{equation}
\begin{equation}
\label{GrindEQ__3_}
R_B \ge \frac{c}{b}r_a + \frac{a}{b}r_c
\end{equation}
\begin{equation}
\label{GrindEQ__4_}
R_C\ge \frac{b}{c}r_a + \frac{a}{c}r_b
\end{equation}
where  $a=\left|BC\right|, \, b=\left|CA\right|, \, c=\left|AB\right|$.

In this paper we determine a set of points \textbf{\textit{E}} for which
\begin{equation}
\label{GrindEQ__5_}
R_A + R_B + R_C
\ge
\left(
\frac{c}{b}+\frac{b}{c}
\right)r_a
+
\left(\
\frac{c}{a}+\frac{a}{c}
\right)
r_b
+
\left(
\frac{a}{b}+\frac{b}{a}
\right)
r_c
\end{equation}
is valid. It is known that the triangular area of $\triangle ABC$ is contained in the set \textbf{\textit{E}}
\cite{[12]}, \cite{[6]}, \cite{[5]}, \cite{[13]}, \cite{[14]}, \cite{[4]}.
Here we show that the set \textbf{\textit{E}} is greater than the triangle
$\triangle ABC$, and we give a geometric interpretation of the set \textbf{\textit{E}}.

The proofs of Erd\" os-Mordell inequality are often based on different proofs of inequality \eqref{GrindEQ__2_}, as given in
\cite{[6]}, \cite{[9]}, \cite{[3]}, \cite{[5]}, \cite{[7]}, \cite{[8]}, \cite{[4]}.
N. Derigades in \cite{[17]} proved the inequality \eqref{GrindEQ__5_} valid in the whole plane of the triangle,
where $r_a, \, r_b$ and $r_c$, are signed distances. A similar result was given by B. Male\v sevi\' c \cite{[10]}, \cite{[11]}.

Note that V. Pambuccian \cite{[21]} recently proved that the Erd\" os-Mordell inequality is equivalent to non-positive curvature.
Overview of recent results on Erd\" os-Mordell inequalities and related inequalities is given in
\cite{[20]} - \cite{[12]}, \cite{[25]}, \cite{[17]}, \cite{[18]}, \cite{[13]} - \cite{[11]}, \cite{[21]}, \cite{[16]}, \cite{[27]} - \cite{[30]} .

\section{The Main Results}

In this section we analyze only the inequality \eqref{GrindEQ__2_}.
Let $\triangle ABC$ be a triangle with vertices $A\left(0,r\right), \, B\left(p,0\right), \, C\left(q,0\right), \, p \ne q, \, r\ne 0$.
Without diminishing generality, let $p<q$. We denote by $M\left(x,y\right)$ an arbitrary point in the plane of the triangle $\triangle ABC$.
The distance from the point $M$ to the point $A$, and the distance from the point $M$ to the straight lines \textbf{\textit{b}}
and \textbf{\textit{c}} are given by functions:
\begin{equation}
\label{GrindEQ__6_}
R_A = \displaystyle\sqrt{\mathstrut x^2+{\left(y-r\right)}^2}
\end{equation}
\begin{equation}
\label{GrindEQ__7_}
r_b = \mbox{\small $\displaystyle\frac{\left|-qy-rx+qr\right|}{\sqrt{\mathstrut r^2+q^2}}$}
\end{equation}
\begin{equation}
\label{GrindEQ__8_}
r_c = \mbox{\small $\displaystyle\frac{\left|py+rx-pr\right|}{\sqrt{r^2+p^2}}$}
\end{equation}
respectively. Consider the inequality \eqref{GrindEQ__2_} related to the vertex $A$. The analytical notation of this inequality is:
\begin{equation}
\label{GrindEQ__9_}
\sqrt{x^2+{\left(y-r\right)}^2}
\ge
\mbox{\small $\displaystyle\frac{\sqrt{r^2+p^2}}{\left|q-p\right|}\frac{\left|-qy-rx+qr\right|}{\sqrt{r^2+q^2}}$}
+
\mbox{\small $\displaystyle\frac{\sqrt{r^2+q^2}}{\left|q-p\right|}\frac{\left|py+rx-pr\right|}{\sqrt{r^2+p^2}}$},
\end{equation}

\break

\noindent
i.e.

\begin{equation}
\label{GrindEQ__10_}
\begin{array}{rcc}
\left|q-p\right| \displaystyle\sqrt{\mathstrut r^2+p^2} \displaystyle\sqrt{\mathstrut r^2+q^2} \displaystyle\sqrt{\mathstrut x^2+{\left(y-r\right)}^2}
&\;\ge&
\left(r^2+p^2\right)\left|-qy-rx+qr\right| \\[1.0 ex]
&\;   &
+
\left(r^2+q^2\right)\left|py+rx-pr\right|.
\end{array}
\end{equation}

Let  $y=kx+r, \, k \in \overline{\mathbb R}$, then the inequality \eqref{GrindEQ__10_} reads as follows:
\begin{equation}
\label{GrindEQ__11_}
\left|x\right|\left|q\!-\!p\right|
\displaystyle\sqrt{\mathstrut r^2\!+\!p^2}
\displaystyle\sqrt{\mathstrut r^2\!+\!q^2}
\displaystyle\sqrt{\mathstrut 1\!+\!k^2}
\ge
\left|x\right|
{\Big (}\!
\left(r^2\!+\!p^2\right)\!\left|-\!qk\!-\!r\right|
+
\left(r^2\!+\!q^2\right)\!\left|pk\!+\!r\right|
\!{\Big )}
\end{equation}

For $x = 0$, the previous inequality is reduced to an equality which solution is the point $A\left(0,r\right)$.
For $x \ne 0$ we obtain inequality by a single variable $k$:
\begin{equation}
\label{GrindEQ__12_}
\left|q\!-\!p\right|
\displaystyle\sqrt{\mathstrut r^2\!+\!p^2}
\displaystyle\sqrt{\mathstrut r^2\!+\!q^2}
\displaystyle\sqrt{\mathstrut 1\!+\!k^2}
\ge
\left(r^2+p^2\right)\!\left|\!-qk\!-\!r\right|
+
\left(r^2+q^2\right)\!\left|pk\!+\!r\right|.
\end{equation}
Solution of the inequality \eqref{GrindEQ__12_} reduces to four cases per parameter $k$:
\begin{equation}
\label{GrindEQ__13_}
\left(\alpha_1\right) \, : \;
\left\{\!
\begin{array}{c}
\;\; pk+r\ge 0 \\
    -qk-r\ge 0,
\end{array}
\right.
\end{equation}
\begin{equation}
\label{GrindEQ__14_}
\left(\alpha_2\right) \, : \;
\left\{\!
\begin{array}{c}
\;\; pk+r<0 \\
    -qk-r\ge 0,
\end{array}
\right.
\end{equation}
\begin{equation}
\label{GrindEQ__15_}
\left(\alpha_3\right) \, : \;
\left\{\!
\begin{array}{c}
\;\; pk+r\ge 0 \\
    -qk-r<0,
\end{array}
\right.
\end{equation}
\begin{equation}
\label{GrindEQ__16_}
\left(\alpha_4\right) \, : \;
\left\{\!
\begin{array}{c}
\;\; pk+r<0 \\
    -qk-r<0.
\end{array}
\right.
\end{equation}

Note that the value $k$ corresponds to the points $(x,y) \in {\mathbb R}^2$ located on the straight line $y=kx+r$.
With its values, the mentioned parameter of the line $y = kx + r$ decomposes ${\mathbb R}^2$ on four corner areas.
Inquiring the existence of parameter $k$ (i.e. the pencil of lines $y=kx+r$ through the vertex $A$) depending on
the signs of parameters $p, \, q$ and $r$, we provide the following table of existing corner areas
$\left(\alpha_1\right)-\left(\alpha_4\right)$:

\begin{center}

\begin{tabular}{|p{0.15in}|p{0.35in}|p{0.35in}|p{0.35in}|p{0.35in}|p{0.35in}|p{0.35in}|p{0.35in}|} \hline
          & $\;\;\;\;\,p$ & $\;\;\;\;\,q$ & $\;\;\;\;\,r$ & $\;\;(\alpha_1)$ & $\;\;(\alpha_2)$ & $\;\;(\alpha_3)$ & $\;\;(\alpha_4)$ \\ \hline
$ \,\;1.$ & $\;\;>$0      & $\;\;>$0      & $\;\;>$0    & \quad\texttt{+}  & \quad\texttt{+}  & \quad\texttt{+}  & \quad\texttt{-}    \\ \hline
$ \,\;2.$ & $\;\;<$0      & $\;\;>$0      & $\;\;>$0    & \quad\texttt{+}  & \quad\texttt{-}  & \quad\texttt{+}  & \quad\texttt{+}    \\ \hline
$ \,\;3.$ & $\;\;<$0      & $\;\;<$0      & $\;\;>$0    & \quad\texttt{-}  & \quad\texttt{+}  & \quad\texttt{+}  & \quad\texttt{+}    \\ \hline
$ \,\;4.$ & $\;\;>$0      & $\;\;>$0      & $\;\;<$0    & \quad\texttt{-}  & \quad\texttt{+}  & \quad\texttt{+}  & \quad\texttt{+}    \\ \hline
$ \,\;5.$ & $\;\;<$0      & $\;\;>$0      & $\;\;<$0    & \quad\texttt{+}  & \quad\texttt{+}  & \quad\texttt{-}  & \quad\texttt{+}    \\ \hline
$ \,\;6.$ & $\;\;<$0      & $\;\;<$0      & $\;\;<$0    & \quad\texttt{+}  & \quad\texttt{+}  & \quad\texttt{+}  & \quad\texttt{-}    \\ \hline
$ \,\;7.$ & $\;\;=0$      & $\;\;>$0      & $\;\;>$0    & \quad\texttt{+}  & \quad\texttt{-}  & \quad\texttt{+}  & \quad\texttt{-}    \\ \hline
$ \,\;8.$ & $\;\;=0$      & $\;\;>$0      & $\;\;<$0    & \quad\texttt{-}  & \quad\texttt{+}  & \quad\texttt{-}  & \quad\texttt{+}    \\ \hline
$ \,\;9.$ & $\;\;<0$      & $\;\;=0$      & $\;\;>$0    & \quad\texttt{-}  & \quad\texttt{-}  & \quad\texttt{+}  & \quad\texttt{+}    \\ \hline
$  10.$   & $\;\;<0$      & $\;\;=0$      & $\;\;<$0    & \quad\texttt{+}  & \quad\texttt{+}  & \quad\texttt{-}  & \quad\texttt{-}    \\ \hline
\end{tabular}

\bigskip

\textit{Table 1: The existence of the corner area depending on the parameters p, q and r}

\end{center}

\break

The corner areas $\left(\alpha_1\right)$ and $\left(\alpha_4\right)$ are always in the interior of
$\sphericalangle BAC$ and its cross angle, while the areas $\left(\alpha_2\right)$ and $\left(\alpha_3\right)$
are in the interior of its supplementary angle (Figure 2).

\begin{center} 

\vspace*{50.0mm} \hspace*{-57.5mm} \includegraphics*[height=50.0mm,keepaspectratio=true,scale=0.80]{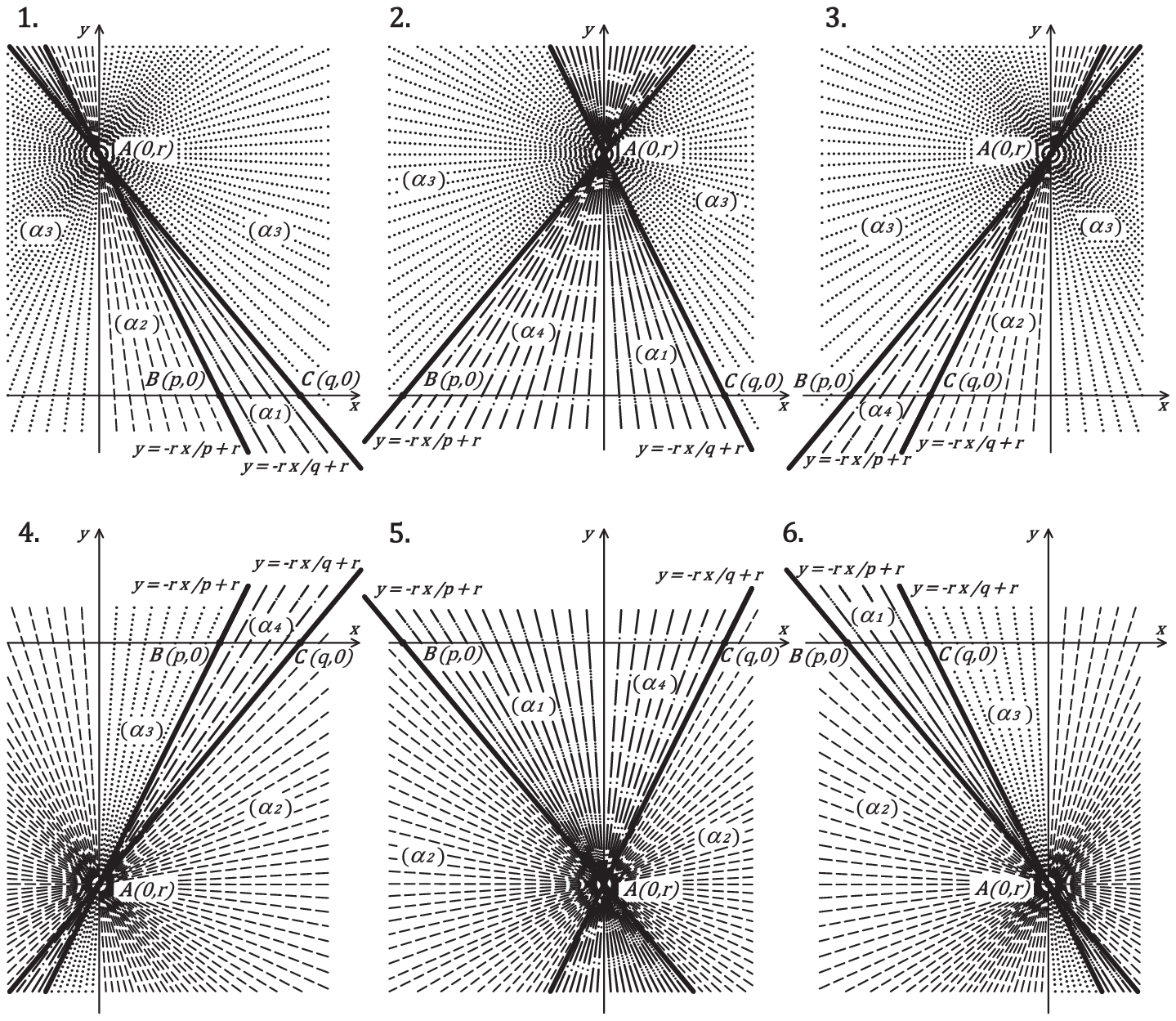}

\noindent
\textit{Figure 2: Existence of the corner area for the vertex A (Cases 1. to 6. in the Table 1)}

\end{center}

\medskip
Let us consider the equation:

\vspace*{-2.0 mm}

\begin{equation}
\label{GrindEQ__17_}
\left(q-p\right)
\displaystyle\sqrt{\mathstrut r^2\!+\!p^2}
\displaystyle\sqrt{\mathstrut r^2\!+\!q^2}
\displaystyle\sqrt{\mathstrut 1\!+\!k^2}
=
\left(r^2\!+\!p^2\right)\left|-qk-r\right|+\left(r^2\!+\!q^2\right)\left|pk+r\right|.
\end{equation}

\medskip 
\noindent
{\large \textbf{I)}}$\;$
Let $k$ fulfill $\left(\alpha_1\right)$ or $\left(\alpha_4\right)$. Then the previous equation can be rewritten in a way that follows,
with positive sign (\texttt{+}) in the case of area $\left(\alpha_1\right)$ and negative sign (\texttt{-}) in the case of area~$\left(\alpha_4\right)$
\begin{equation}
\label{GrindEQ__18_}
\left(q-p\right)
\displaystyle\sqrt{\mathstrut r^2\!+\!p^2}
\displaystyle\sqrt{\mathstrut r^2\!+\!q^2}
\displaystyle\sqrt{\mathstrut 1\!+\!k^2}
=
\pm
{\big (}\!\left(-qk-r\right)\left(r^2\!+\!p^2\right)+\left(pk+r\right)\left(r^2\!+\!q^2\right)\!{\big )}
\end{equation}
i.e.
\begin{equation}
\label{GrindEQ__19_}
\left(q-p\right)
\displaystyle\sqrt{\mathstrut r^2\!+\!p^2}
\displaystyle\sqrt{\mathstrut r^2\!+\!q^2}
\displaystyle\sqrt{\mathstrut 1\!+\!k^2}
=
\pm
\left(q-p\right){\big (}r\left(q\!+\!p\right)+k\left(pq-r^2\right)\!{\big )}
\end{equation}
abbreviated written as

\break

\noindent
\begin{equation}
\label{GrindEQ__20_}
\lambda
\displaystyle\sqrt{\mathstrut 1+k^2}
=
\pm \beta k \pm \gamma
=\left\{
\begin{array}{c}
\;\;\;\;\beta k+\gamma ,\ \ \ \ k \in \left(\alpha_1\right) \\[1.0 ex]
       -\beta k-\gamma ,\ \ \   k \in \left(\alpha_4\right)
\end{array}
\right.
\end{equation}
where at:
\begin{equation}
\label{GrindEQ__21_}
\lambda
=
\left(q-p\right)
\displaystyle\sqrt{\mathstrut r^2+p^2}
\displaystyle\sqrt{\mathstrut r^2+q^2}
\;\;\;\;\mbox{\rm and}\;\;\;\;\lambda>0
\end{equation}
\begin{equation}
\label{GrindEQ__22_}
\beta
=
\left(pq-r^2\right)\left(q-p\right)
\end{equation}
\begin{equation}
\label{GrindEQ__23_}
\gamma
=
r\left(q^2-p^2\right).
\end{equation}
As $p<q$, the equation \eqref{GrindEQ__19_} can be divided by $ q - p \ne 0$ and then squared:
\begin{equation}
\label{GrindEQ__24_}
\left(r^2\!+\!p^2\right)\left(r^2\!+\!q^2\right)\left(1+k^2\right)
=
{\big (} \!\, r \left(q+p\right) + k \left( pq-r^2 \right) \! {\big )}^2
\end{equation}
which transforms into
\begin{equation} \label{GrindEQ__25_}
{\big (} \!\, r\left(p+q\right)k-\left(pq-r^2\right) \! {\big )}^2 = 0.
\end{equation}
Based on the above equation, we conclude that there exists the unique solution:
\begin{equation} \label{GrindEQ__26_}
k_1=\frac{pq-r^2}{r\left(p+q\right)}
\end{equation}
only if, for $k=k_1$:
\begin{equation} \label{GrindEQ__27_}
\pm \beta k \pm \gamma \ge 0
\end{equation}
is valid.

\smallskip
Hence, the straight line $y=k_1x+r$ is in the interior of $\sphericalangle BAC$ and its cross angle, or it doesn't exist.
The cases where values $k_1$ from the formula \eqref{GrindEQ__26_} does not meet the condition \eqref{GrindEQ__27_}
are presented in the \textit{Table 1 }with:

\smallskip
in the case 1: $k_{{\rm 1}}{\rm >}-{r}/{q}\Longleftrightarrow p\left(q^{{\rm 2}}{\rm +}r^{{\rm 2}}\right){\rm >}0\ $;

\smallskip
in the case 3: $k_{{\rm 1}}{\rm >}-{r}/{p}\Longleftrightarrow \left({\rm -}q\right)\left(p^{{\rm 2}}{\rm +}r^{{\rm 2}}\right){\rm >}0$;

\smallskip
in the case 4: $k_{{\rm 1}}{\rm <}-{r}/{q}\Longleftrightarrow p\left(q^{{\rm 2}}{\rm +}r^{{\rm 2}}\right){\rm >}0\ $;

\smallskip
in the case 6: $k_{{\rm 1}}{\rm <}-{r}/{p}\Longleftrightarrow \left({\rm -}q\right)\left(p^{{\rm 2}}{\rm +}r^{{\rm 2}}\right){\rm >}0$.

\begin{lemma}
For $k\in \left(\alpha_1\right)\ \cup \left(\alpha_4\right)$ inequality \eqref{GrindEQ__12_} is valid,
where equality holds for $k=k_1$ if \eqref{GrindEQ__27_} is fulfilled.
\end{lemma}

\vspace*{-3.0 mm}

\begin{proof}
$\mbox{\eqref{GrindEQ__12_}} \Longleftrightarrow {\big (} r\left(p+q\right)k - \left(pq-r^2\right){\big )}^2 \ge 0$. $\Box$
\end{proof}

\begin{corollary}
Inequality \eqref{GrindEQ__12_} is valid for lines \textbf{\textit{b}} and \textbf{\textit{c}}.
\end{corollary}

\medskip 
\noindent
{\large \textbf{II)}}$\;$
Let $k$ fulfill $\left(\alpha_2\right)$ or $\left(\alpha_3\right)$. Then equation \eqref{GrindEQ__17_} can be rewritten
in a way that follows, with negative sign (\texttt{-}) in the case of area $\left(\alpha_2\right)$ and positive sign (\texttt{+})
in the case of area $\left(\alpha_3\right)$
\begin{equation}
\label{GrindEQ__28_}
\left(q-p\right)
\displaystyle\sqrt{\mathstrut r^2\!+\!p^2}
\displaystyle\sqrt{\mathstrut r^2\!+\!q^2}
\displaystyle\sqrt{\mathstrut 1\!+\!k^2}
=
\pm
{\big (}\!\left(qk+r\right)\left(r^2\!+\!p^2\right) + \left(pk+r\right)\left(r^2\!+\!q^2\right)\!{\big )}
\end{equation}
or abbreviated written as
\begin{equation}
\label{GrindEQ__29_}
\lambda
\displaystyle\sqrt{\mathstrut 1 + k^2}
=
\pm \delta k \pm \varepsilon
=
\left\{
\begin{array}{c}
\;\;\; \delta k + \varepsilon,\ \ \ \ k\in \left(\alpha_3\right) \\[1.0 ex]
      -\delta k - \varepsilon,\ \ \   k\in \left(\alpha_2\right)
\end{array}
\right.
\end{equation}
with parameters:
\begin{equation*}
\lambda
=
\left(q-p\right)
\displaystyle\sqrt{\mathstrut r^2+p^2}
\displaystyle\sqrt{\mathstrut r^2+q^2} \quad \mbox{and} \quad \lambda >0
\end{equation*}
\begin{equation}
\label{GrindEQ__30_}
\delta
=
\left(r^2+pq\right)\left(q+p\right)
\end{equation}
\begin{equation}
\label{GrindEQ__31_}
\varepsilon
=
r\left(2r^2+q^2+p^2\right) \! .
\end{equation}

The equation \eqref{GrindEQ__29_} is considered under the following condition:
\begin{equation}
\label{GrindEQ__32_}
\pm \delta k\pm \varepsilon
\ge
0.
\end{equation}

By squaring the equation \eqref{GrindEQ__29_} we obtain
\begin{equation}
\label{GrindEQ__33_}
P(k)
=
\lambda^2\left(1+k^2\right) - \left(\pm \delta k \pm \varepsilon \right)^2
=
\left(\lambda^2-\delta^2\right) k^2 - 2 \delta \varepsilon k + \left(\lambda^2-\varepsilon^2\right)
=0.
\end{equation}
For the square trinomial
\begin{equation}
\label{GrindEQ__34_}
P(k)
=
\widehat{{\rm A}} \, k^2 + \widehat{{\rm B}} \, k + \widehat{{\rm C}}
\end{equation}
coefficients $\widehat{{\rm A}}, \, \widehat{{\rm B}}, \, \widehat{{\rm C}}$ are determined by:
\begin{equation}
\label{GrindEQ__35_}
\widehat{{\rm A}}
=
\lambda^2 - \delta^2
=
\left(q-p\right)^2 \! \left(r^2\!+\!p^2\right) \! \left(r^2\!+\!q^2\right) - \left(r^2+pq\right)^2 \! \left(q+p\right)^2
\end{equation}
\begin{equation}
\label{GrindEQ__36_}
\widehat{{\rm B}}
=
-2 \delta \varepsilon
=
-2 r\left(r^2 + pq\right) \! {\big (}q + p{\big )} \! \left(2 r^2 + q^2 + p^2 \right)
\end{equation}
\begin{equation}
\label{GrindEQ__37_}
\widehat{{\rm C}}
=
\lambda^2 - \varepsilon^2
=
\left(r^2+pq\right){\big (}\!\left(pq-r^2\right) \left(q-p\right)^2 - 2r^2(2r^2+q^2+p^2){\big )}.
\end{equation}
Let us consider the equation:
\begin{equation}
\label{GrindEQ__38_}
\widehat{{\rm A}}
=
-4pqr^4 + \left(p^4\!+\!q^4\!-\!4pq^3\!-\!4p^3q\!-\!2p^2q^2\right) \! r^2 - 4p^3q^3
=
0.
\end{equation}
It has real solutions for  $r$ in the following form:
\begin{equation}
\label{GrindEQ__39_}
\left\{
\begin{array}{c}
r_{1,2}
=
\mbox{\small \mbox{$\displaystyle\frac{1}{4\displaystyle\sqrt{\mathstrut pq}}$}}
\left((q-p)^2\ \pm \ \displaystyle\sqrt{\mathstrut {\left(q-p\right)}^4-16{p^2q}^2} \;\; \right) > 0 \\[1.75 ex]
r_{3,4}
=
-\mbox{\small \mbox{$\displaystyle\frac{1}{4\displaystyle\sqrt{\mathstrut pq}}$}}
\left((q-p)^2\ \pm \ \displaystyle\sqrt{\mathstrut {\left(q-p\right)}^4-16{p^2q}^2} \;\; \right) < 0
\end{array}
\right.
\end{equation}
iff
\begin{equation}
\label{GrindEQ__40_}
{\Big (} p \ge 0 \; \wedge \; q \ge (3+2\displaystyle\sqrt{2})p {\Big )}
\;\; \vee \;\;
{\Big (} p < 0 \; \wedge \; q \le (3-2\displaystyle\sqrt{2})p {\Big )}.
\end{equation}

\begin{remark}
When $p < 0$ and $q > 0$ then $\widehat{{\rm A}}
=
4\left|p\right|qr^4 + \left(q^2{-p}^2\right)^2 \! r^2 + 4\left|p\right|q$ $\left(p^2+q^2\right) \! r^2 + 4\left|p\right|^3q^3 > 0$ is valid.
Note that the equation $\widehat{{\rm A}} =0$ is not considered for $p=0$ or $q=0$ {\big (}because we obtain the contradictions:
$p = 0, \; q \ne 0$:  $\widehat{{\rm A}} = r^2 q^4 = 0 \, \Longrightarrow \,  r=0$; i.e. $p \ne 0, \; q=0$:
$\widehat{{\rm A}} = r^2 p^4 = 0 \, \Longrightarrow \, r = 0${\big )}.
\end{remark}

\noindent

We distinguish the cases:

\medskip 
\noindent
{\large \textbf{a)}}
Let $r = r_j$ for some $j = 1, \, 2, \, 3, \, 4$, then $\widehat{{\rm A}} = 0$. In this case, $\widehat{{\rm B}} \ne 0$,
because $r^2 + pq \ne 0$ and $q + p \ne 0$ {\big (}in the case of equilateral triangle, there will be valid $q + p = 0$
and then $r = \pm pi, \; i = \displaystyle\sqrt{\!-\!1}\,${\big )}. Therefore, by solving the linear equation
$\widehat{{\rm B}} \, k + \widehat{{\rm C}} = 0$ we find that:
\begin{equation}
\label{GrindEQ__41_}
k_2
=
-\displaystyle\frac{\widehat{{\rm C}}}{\widehat{{\rm B}}}
=
\displaystyle\frac{\lambda^2 - \varepsilon^2}{2 \delta \varepsilon}
=
\displaystyle\frac{
\left(q-p\right)^2\left(r^2+p^2\right)\left(r^2+q^2\right) - r^2\left(2r^2+q^2+p^2\right)^2
}{
2r\left(q+p\right)\left(2r^2+q^2+p^2\right)
}.
\end{equation}

For $p < 0$ and $q > 0$ the case \textbf{a)} is not considered {\big (}because $\widehat{{\rm A}} > 0${\big )}. Let us examine when the value
$k_2$ meet the condition \eqref{GrindEQ__32_}. It is valid that:
\begin{equation*}
\pm \delta k_2 \pm \varepsilon
\ge
0
\; \Longleftrightarrow \;
\pm \left(\delta k_2+\varepsilon \right)
=
\pm \left(\delta \displaystyle\frac{\,\lambda^2 - \varepsilon^2}{2 \delta \varepsilon} + \varepsilon \right)
=
\pm \left(\displaystyle\frac{\,\lambda^2+\varepsilon^2}{2 \varepsilon}\right)
\ge
0.
\end{equation*}
Based on $\varepsilon = r\left(2r^2+q^2+p^2\right)$ we conclude:

\smallskip
\noindent
if $r > 0$ then $\delta k_2 + \varepsilon \ge 0$ is fulfilled, whereby $k_2$ fulfills condition \eqref{GrindEQ__32_}
and $k_2 \in \left(\alpha_3\right)$;

\smallskip
\noindent
if $r < 0$ then $-\delta k_2-\varepsilon \ge 0\ $is fulfilled, whereby $k_2$ fulfills condition \eqref{GrindEQ__32_}
and $k_2 \in \left(\alpha_2\right)$.

\smallskip
\noindent
In this case, the line $y = k_2 x + r$ is in the exterior of $\sphericalangle BAC$ and its cross angle.

\medskip 
\noindent
{\large \textbf{b)}}
Let $r \ne r_j$ for each $j = 1, \, 2, \, 3, \, 4$, then $\widehat{{\rm A}} \ne 0$ and in this case, by solving the quadratic equation
\eqref{GrindEQ__33_}, we find the values:
\begin{equation}
\label{GrindEQ__42_}
\begin{array}{rcl}
k_{2,3}
& \!\!=\!\! &
\displaystyle\frac{
-\delta \varepsilon
\pm \displaystyle\sqrt{\mathstrut \lambda^2 \left(\delta^2 \!+\! \varepsilon^2 \!-\! \lambda^2\right)}}{\delta^2\!-\!\lambda^2}  \\[1.0 ex]
& \!\!=\!\! &
\displaystyle\frac{r(p+q)(r^2\!+pq)(q^2\!+\!p^2\!+\!2r^2)
\pm 2\left(r^2\!+\!p^2\right)\left(r^2\!+\!q^2\right)\left(q\!-\!p\right)\displaystyle\sqrt{\mathstrut r^2\!+pq}}{
\left(q\!-\!p\right)^2\left(r^2\!+\!p^2\right)\left(r^2\!+\!q^2\right)-\left(r^2\!+pq\right)^2\left(q\!+\!p\right)^2}.
\end{array}
\end{equation}
If $r^2 + pq \ge 0$ then exists $k_{2,3}\in {\mathbb R}$. Incidence of $k_{2,3} \in {\mathbb R}$ to the area
$\left(\alpha_3\right)$, as to the area $\left(\alpha_2\right)$ is determined by the inequality \eqref{GrindEQ__32_}.
The expression $\delta k_{2,3} + \varepsilon$ exists for $\delta \ne \pm \lambda $, whereby the expression
$\delta k_{2,3} + \varepsilon $ is either positive or negative (because $\delta k_{2,3} + \varepsilon = 0 \Longrightarrow \delta = \pm \lambda $).

\smallskip
\noindent
Based on the Corollary 1, the straight lines $y=k_sx+r, (s=2,\ 3)$ are in the exterior of $\sphericalangle BAC$ and its cross angle (\textit{Figure 3}).

\smallskip
\noindent
Consider the limiting case for $k_{2,3}$ when $r \to r_{j}$. Note that $\widehat{{\rm A}} = \lambda^2 - \delta^2
{{\underset{r\to r_{j}}{\longrightarrow}}} \,0$ is valid, whereat from
\begin{equation*}
k_{2,3}
=
\displaystyle\frac{-\varepsilon}{\left(\delta -\lambda\right)\left(\delta +\lambda\right)}
\cdot
\left(
\delta
\mp
\left|\lambda\right|
\displaystyle\sqrt{\mathstrut 1 + \frac{\delta^2-\lambda^2}{\varepsilon^2}}\;\;
\right)
\end{equation*}
follows
\begin{equation*}
\mathop{\mbox{\rm lim}}_{r\to r_j}{k}_2
=
\displaystyle\frac{-\varepsilon }{\left(\delta + \lambda\right)}
\;\;\; \wedge \;\;\;
\mathop{\mbox{\rm lim}}_{r\to r_j}{k}_3 = \infty .
\end{equation*}

\break

\noindent
\begin{center} 

\vspace*{30.0mm} \hspace*{-47.5mm} \includegraphics*[height=60.0mm,keepaspectratio=true,scale=0.54]{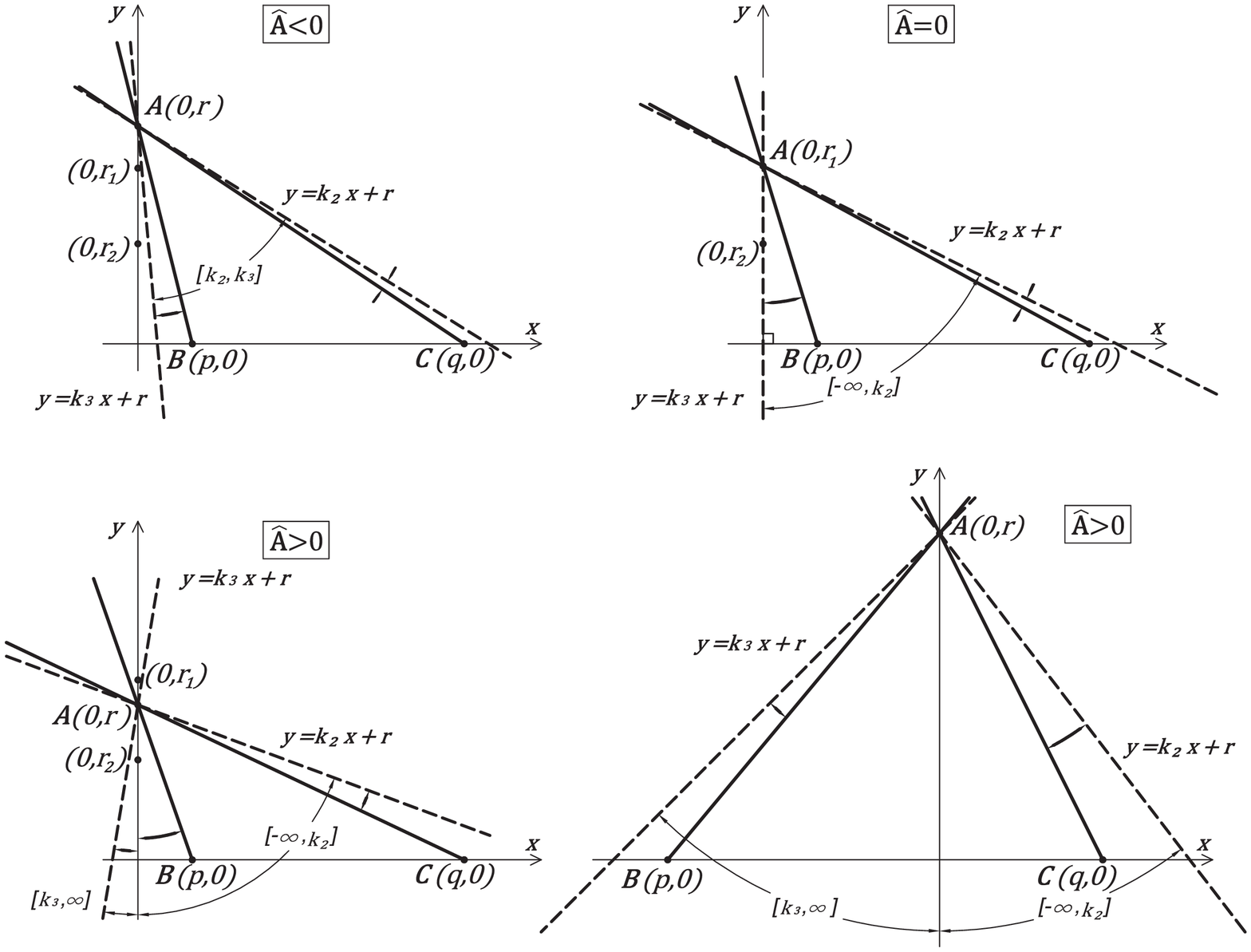}

\noindent
\textit{Figure 3: The existence of the lines $y = k_s x \!+\! r,\; (\!\mbox{\small $s\!=\!2,3$}\!)$ depending~on~the~parameter~$\widehat{{\rm A}}$}

\end{center}

\smallskip
\noindent
Related to the $\sphericalangle BAC$ we distinguish the cases:

\bigskip
\noindent
\textbf{1.}
$\sphericalangle BAC < \pi/2 \Longleftrightarrow r^2 + pq > 0$ and if $\widehat{{\rm A}} \ne 0$ then there are two real and different values of
$k_2$ and $k_3$. In this case, the following lemma is valid:

\begin{lemma}
For $\sphericalangle BAC < \pi/2$, $k \in \left(\alpha_2\right) \cup \left(\alpha_3\right)$ the inequality
\eqref{GrindEQ__12_} is valid, just in the cases:

\begin{enumerate}
\item
$
\widehat{{\rm A}} > 0
\; \wedge \;
k \in \left[-\infty, \, k_2\right] \cup \left[k_3, \, +\infty\right] \, \backslash \, {\big (}\!\left(\alpha_1\right) \cup \left(\alpha_4\right)\!{\big )}
$;

\item
$
\widehat{{\rm A}} = 0
\; \wedge \;
k \in \left[-\infty, \, k_2\right] \backslash {\big (}\left(\alpha_1\right) \cup \left(\alpha_4\right){\big )}
$;

\item
$
\widehat{{\rm A}} < 0
\; \wedge \;
k \in \left[k_2, \, k_3\right] \backslash {\big (}\left(\alpha_1\right) \cup \left(\alpha_4\right){\big )}
$;
\end{enumerate}

\noindent
where the equality holds for $k = k_2$ or $k = k_3$.
\end{lemma}

\bigskip
\noindent
\textbf{2.}
If $\sphericalangle BAC = \pi/2 \Longleftrightarrow r^2+pq=0$ then $\widehat{{\rm A}} = -qp\left(q-p\right)^4$,
$\widehat{{\rm B}} = 0 $ and $\widehat{{\rm C}} = 0$, according to the equation \eqref{GrindEQ__42_} that $k_{2,3}=0.$
Hence is valid:

\begin{lemma}
For $\sphericalangle BAC=\pi/2$ and $k \in \left(\alpha_2\right) \cup \left(\alpha_3\right)$ the inequality \eqref{GrindEQ__12_} is valid.
The equality is valid only for $k = 0$.
\end{lemma}

\vspace*{-3.0 mm}

\begin{proof}
$\mbox{\eqref{GrindEQ__12_}}
\Longleftrightarrow
\widehat{{\rm A}} \, k^2 + \widehat{{\rm B}} \, k + \widehat{{\rm C}} \ge 0
\Longleftrightarrow
-qp\left(q-p\right)^4 k^2 \ge 0$. $\Box$
\end{proof}

\break

\bigskip
\noindent
\textbf{3.}
$\sphericalangle BAC > \pi/2 \Longleftrightarrow r^2+pq < 0$. In this case, for: $r^2 < -pq$ and for the coefficient~$\widehat{{\rm A}}$:
\begin{equation*}
\begin{array}{rcl}
\widehat{{\rm A}}
& \!>\! &
4r^6 + \left(p^4+q^4\right)r^2 + 4\left(p^2+q^2\right)r^4 - 2r^6 + 4p^2q^2r^2                     \\[1.25 ex]
& \!=\! &
2r^6 + 4\left(p^2+q^2\right)r^4 + \left(p^4+q^4+4p^2q^2\right)r^2
>
0
\end{array}
\end{equation*}
is valid. Since $k_{2,3} \in \mathbb C$ and $\widehat{{\rm A}} > 0$ the inequality \eqref{GrindEQ__12_} is valid, which proves the claim:

\vspace*{-0.5 mm}

\begin{lemma}
For $\sphericalangle BAC > \pi/2$ and $k \in \left(\alpha_2\right) \cup \left(\alpha_3\right)$ the inequality \eqref{GrindEQ__12_}
is valid in the strict form.
\end{lemma}

Based on the previous considerations in \textbf{I)} and \textbf{II)}, follows:

\vspace*{-0.5 mm}

\begin{statement}
The inequality \eqref{GrindEQ__12_} holds in following cases:
\begin{equation*}
\begin{array}{lcl}
                 & \hspace*{20 mm} & k \in \left(\alpha_1\right) \cup \left(\alpha_4\right)                        \\[1.0 ex]
\hspace*{-15.75 mm}
\mbox{\it or}    & \hspace*{20 mm} &                                                                               \\[1.0 ex]
                 & \hspace*{20 mm} & k \in \left(\alpha_2\right) \cup \left(\alpha_3\right)
                                     \;\;\mbox{\it for}\;\; \sphericalangle BAC \ge \pi/2                          \\[1.0 ex]
\hspace*{-15.75 mm}
\mbox{\it i.e.}  & \hspace*{20 mm} &                                                                               \\[1.0 ex]
                 & \hspace*{20 mm} & k \in \left[-\infty, \, k_2\right] \cup \left[k_3, \, +\infty\right]
                                     \backslash {\big (}\left(\alpha_1\right) \cup \left(\alpha_4\right){\big )}
                                     \;\wedge\; \widehat{{\rm A}} > 0                                              \\[1.0 ex]
                 & \hspace*{20 mm} & k \in \left[-\infty, \, k_2\right]
                                     \backslash {\big (}\left(\alpha_1\right) \cup \left(\alpha_4\right){\big )}
                                     \;\wedge\; \widehat{{\rm A}} = 0                                              \\[1.0 ex]
                 & \hspace*{20 mm} & k \in \left[k_2, \, k_3\right]
                                     \backslash {\big (}\left(\alpha_1\right) \cup \left(\alpha_4\right){\big )}
                                     \;\wedge\; \widehat{{\rm A}} <0,
\end{array}
\end{equation*}
for $\sphericalangle BAC  < \pi/2$.
\end{statement}

\section{Conclusion}

For the vertex $A$, let us define
\begin{equation*}
\mbox{\textbf{\textit{E}}}_A
=
\left\{
\left( x, y \right) \mbox{\large $|$}
R_A
\ge
\mbox{\small $\displaystyle\frac{c}{a}$} r_b + \mbox{\small $\displaystyle\frac{b}{a}$} r_c
\right\},
\end{equation*}
and for the vertices $B$ and $C$, let us define
\begin{equation*}
\mbox{\textbf{\textit{E}}}_B
=
\left\{
\left( x, y \right) \mbox{\large $|$}
R_B
\ge
\mbox{\small $\displaystyle\frac{c}{b}$} r_a + \mbox{\small $\displaystyle\frac{a}{b}$} r_c
\right\},
\end{equation*}

\vspace*{-2.0 mm}

\begin{equation*}
\mbox{\textbf{\textit{E}}}_C
=
\left\{
\left( x, y \right) \mbox{\large $|$}
R_C
\ge
\mbox{\small $\displaystyle\frac{b}{c}$} r_a + \mbox{\small $\displaystyle\frac{a}{c}$} r_b
\right\},
\end{equation*}
respectively. Based on the analysis of the inequalities \eqref{GrindEQ__2_}, \eqref{GrindEQ__3_}
and \eqref{GrindEQ__4_}, the inequality \eqref{GrindEQ__5_} is valid in the intersection of the areas:
\begin{equation}
\label{GrindEQ__43_}
\mbox{\textbf{\textit{E}}}
=
\mbox{\textbf{\textit{E}}}_A
\cap
\mbox{\textbf{\textit{E}}}_B
\cap
\mbox{\textbf{\textit{E}}}_C.
\end{equation}

\vspace*{-1.5 mm}

\noindent
Therefore follows

\vspace*{-0.5 mm}

\begin{statement}
Erd\" os-Mordell inequality is valid in the area $\mbox{\textbf{\textit{E}}}\,$.
\end{statement}

\vspace*{-0.5 mm}

Let us define the set $\textbf{\textit{M}}$ by the intersection of the corner areas formed from
$\mbox{\textbf{\textit{E}}}_A$, $\mbox{\textbf{\textit{E}}}_B$ and $\mbox{\textbf{\textit{E}}}_C$,
containing the initial triangle. Then the set of points $\textbf{\textit{M}}$ is quadrilateral or hexagonal shape,
and is contained the area $\mbox{\textbf{\textit{E}}}$ (\textit{Figure 4}).

\break

\noindent
\begin{center} 

\vspace*{4mm} \hspace*{-35.0mm} \includegraphics*[height=25.0mm,keepaspectratio=true,scale=0.70]{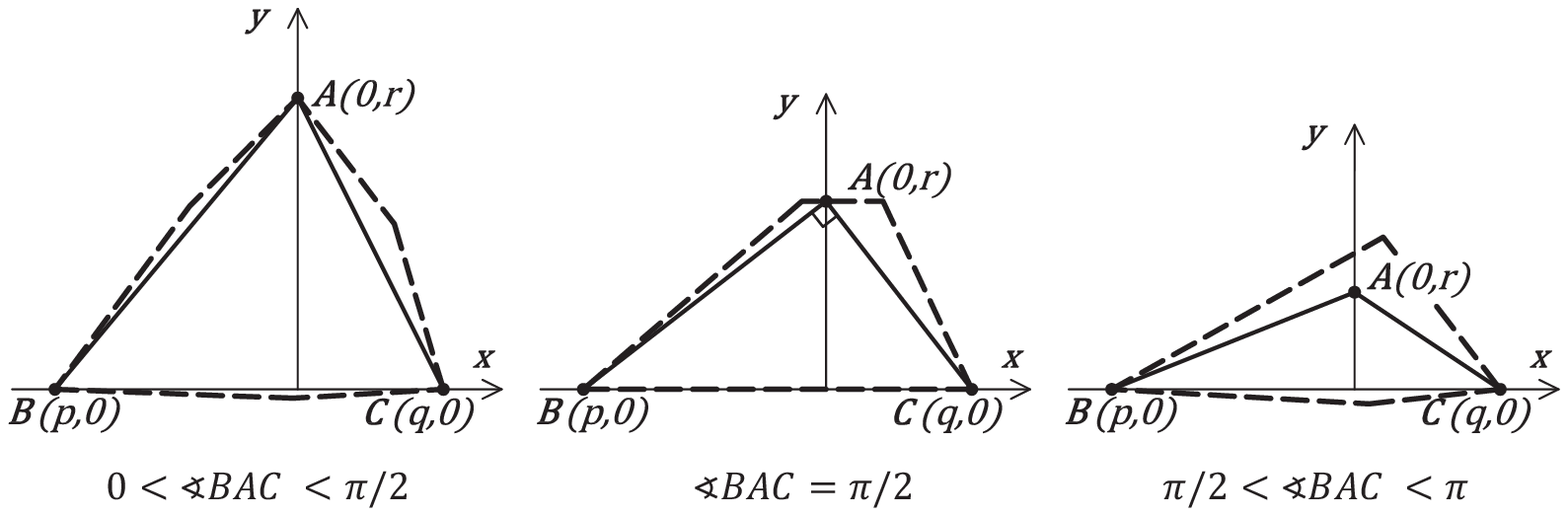}

\smallskip
\noindent
\textit{Figure 4: Extension of the triangle ABC to the area ${\mathbf M} \subset {\mathbf E}$}

\end{center}

\smallskip
Let us define Erd\" os-Mordell curve in the plane of triangle, by the following equation:
\begin{equation}
\label{GrindEQ__44_}
R_A + R_B + R_C
=
2\left(r_a + r_b + r_c\right),
\end{equation}
where
\begin{equation*}
\begin{array}{ccccc}
R_A
=
\displaystyle\sqrt{\mathstrut x^2+{\left(y-r\right)}^2} \, ,
& \;\;\; &
R_B
=
\displaystyle\sqrt{\mathstrut {\left(x-p\right)}^2+y^2} \, ,
& \;\;\; &
R_C
=
\displaystyle\sqrt{\mathstrut {\left(x-q\right)}^2+y^2} \, ,                                                 \\[1.5 ex]
r_a
=
\mbox{\small $\displaystyle\frac{\left|y\left(q-p\right)\right|}{\displaystyle\sqrt{\mathstrut {\left(q-p\right)}^2}}$}
=
\left|y\right| \, ,
& \;\;\; &
r_b
=
\mbox{\small $\displaystyle\frac{\left|-q\left(y-r\right)-rx\right|}{\displaystyle\sqrt{\mathstrut r^2+q^2}}$} \, ,
& \;\;\; &
r_c
=
\mbox{\small $\displaystyle\frac{\left|-p\left(y-r\right)-rx\right|}{\displaystyle\sqrt{\mathstrut r^2+p^2}}$} \, .
\end{array}
\end{equation*}
The curve \eqref{GrindEQ__44_} is a union of parts of algebraic curves of order eight (Figure 5).

\begin{center} 

\vspace*{35.0mm} \hspace*{-45.0mm} \includegraphics*[height=60.0mm,keepaspectratio=true,scale=0.57]{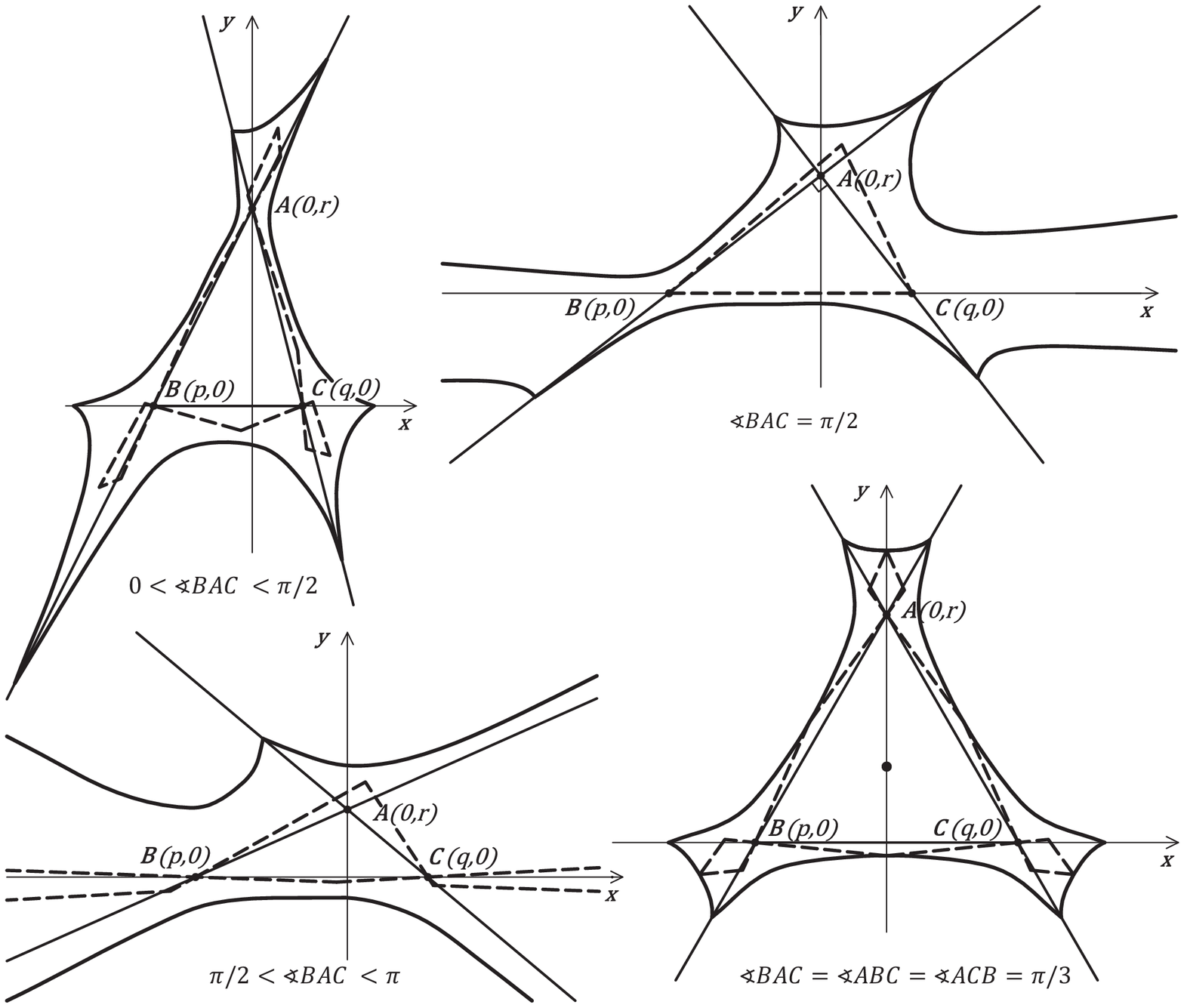}

\smallskip
\noindent
\textit{Figure 5: Erd\" os-Mordell curve and the area ${\mathbf E}$}

\end{center}

\break

Let us denote by \mbox{\textbf{\textit{E'}}} the part of the plane ${\mathbb R}^2$ bounded by
the Erd\" os-Mordell's curve and consisting the triangle $\triangle ABC$. Thus, according to
the fact that inequality \eqref{GrindEQ__5_} is valid in the area of the triangle $\triangle ABC$,
and based on continuity, it follows that inequality \eqref{GrindEQ__5_} is valid in the area \mbox{\textbf{\textit{E'}}}.
Remark that the area \textbf{\textit{E'}} allows us to precise when, except for the inequality \eqref{GrindEQ__5_},
some of the inequalities \eqref{GrindEQ__2_}, \eqref{GrindEQ__3_} and/or \eqref{GrindEQ__4_} are true.
For example, in the area $\left(\mbox{\textbf{\textit{E'}}} \backslash \mbox{\textbf{\textit{E}}}_A\right)
\cap \mbox{\textbf{\textit{E}}}_B \cap \mbox{\textbf{\textit{E}}}_C$ the inequalities \eqref{GrindEQ__5_},
\eqref{GrindEQ__4_}, \eqref{GrindEQ__3_} are true and \eqref{GrindEQ__2_} is not true. At end of this section
let us emphasize that the following statement is true.

\vspace*{-3.0 mm}

\begin{statement}
All geometric inequalities based on the inequalities \eqref{GrindEQ__2_}, \eqref{GrindEQ__3_} and \eqref{GrindEQ__4_}
can be extended from the triangle interior to the area $\mbox{\textbf{\textit{E}}}\,$.
\end{statement}

\vspace*{-5.0 mm}

\begin{example}
In the area $\mbox{\textbf{\textit{E}}}$, the inequality of Child \cite{[3]} is valid:
\begin{equation}
\label{GrindEQ__45_}
R_A \cdot R_B \cdot R_C
\ge
8 \cdot r_a \cdot r_b \cdot r_c
\end{equation}
because, based on inequality between arithmetic and geometric mean, follows:
\begin{equation}
\label{GrindEQ__46_}
a \cdot R_A \ge b \cdot r_c + c \cdot r_b \ge 2\displaystyle\sqrt{\mathstrut b \cdot c \cdot r_b \cdot r_c}
\end{equation}

\vspace*{-2.75 mm}

\begin{equation}
\label{GrindEQ__47_}
b \cdot R_B \ge c \cdot r_a + a \cdot r_c \ge 2\displaystyle\sqrt{\mathstrut c \cdot a \cdot r_c \cdot r_a}
\end{equation}

\vspace*{-2.5 mm}

\begin{equation}
\label{GrindEQ__48_}
c \cdot R_C \ge a \cdot r_b + b \cdot r_a \ge 2\displaystyle\sqrt{\mathstrut a \cdot b \cdot r_a \cdot r_b}.
\end{equation}
Hence, by multiplying the left and right sides of inequalities \eqref{GrindEQ__46_} - \eqref{GrindEQ__48_},
we get the inequality \eqref{GrindEQ__45_} in the area $\mbox{\textbf{\textit{E}}}$. $\Box$
\end{example}

At the end of this paper, let us set up an open problem (proposed by anonymous reviewer):
prove or disprove that there exist a positive number $\varepsilon$ such that the area of
$\mbox{\textbf{\textit{E'}}}$ is bigger than 1+$\varepsilon $ times the area of the triangle for every triangle.
Thus, we set a conjecture: for the finite area of $\mbox{\textbf{\textit{E'}}}$ the value $\varepsilon$
is determined in the case of equilateral triangle.

\bigskip
{\sc Acknowledgment.}
The authors would like to thank anonymous reviewer for his/her valuable comments and suggestions, which were helpful in improving the paper.

\bigskip

\bigskip

\bigskip

\noindent {\footnotesize \em                                                                                      %
\hspace*{5.0mm}  Branko Male\v sevi\' c $($\mbox{\scriptsize \tt CORRESPONDING AUTHOR}$)$         \\[-0.5 ex]     %
\hspace*{5.0mm}  Faculty of Electrical Engineering, University of Belgrade,                       \\[-0.5 ex]     %
\hspace*{5.0mm}  Bulevar Kralja Aleksandra 73, 11000 Belgrade, Serbia                             \\[-0.5 ex]     %
\hspace*{5.0mm}  \email{malesevic@etf.rs}                                                         \\[ 3.0 ex]     %
\hspace*{5.0mm}  Maja Petrovi\' c $($\mbox{\scriptsize \tt CORRESPONDING AUTHOR}$)$               \\[-0.5 ex]     %
\hspace*{5.0mm}  Faculty of Transport and Traffic Engineering, University of Belgrade,            \\[-0.5 ex]     %
\hspace*{5.0mm}  Vojvode Stepe 305, 11000 Belgrade, Serbia                                        \\[-0.5 ex]     %
\hspace*{5.0mm}  \email{majapet@sf.bg.ac.rs}                                                      \\[ 3.0 ex]     %
\hspace*{5.0mm}  Marija Obradovi\' c,                                                             \\[-0.5 ex]     %
\hspace*{5.0mm}  Faculty of Civil Engineering, University of Belgrade,                            \\[-0.5 ex]     %
\hspace*{5.0mm}  Bulevar Kralja Aleksandra 73, 11000 Belgrade, Serbia                             \\[-0.5 ex]     %
\hspace*{5.0mm}  \email{marijao@grf.bg.ac.rs}                                                     \\[ 0.0 ex]     %
\newpage                                                                                                          %
\noindent                                                                                                         %
\hspace*{5.0mm}  Branislav Popkonstantinovi\' c,                                                  \\[-0.5 ex]     %
\hspace*{5.0mm}  Faculty of Mechanical Engineering, University of Belgrade,                       \\[-0.5 ex]     %
\hspace*{5.0mm}  Kraljice Marije 16, 11000 Belgrade, Serbia                                       \\[-0.5 ex]     %
\hspace*{5.0mm}  Faculty of Technical Sciences, University of Novi Sad,                           \\[-0.5 ex]     %
\hspace*{5.0mm}  Trg D. Obradovi\' ca 16, 21000 Novi Sad, Serbia                                  \\[-0.5 ex]     %
\hspace*{5.0mm}  \email{dr.branislav.pop@gmail.com}}                                              \\[ 0.0 ex]     %
                                                                                                                  %


\end{document}